\documentclass[12pt]{amsart}
\usepackage{amsfonts, amsmath}

\usepackage{graphics}
\usepackage{amssymb}
\usepackage{amsmath}
\usepackage{latexsym}
\usepackage{amsthm}
\usepackage[all]{xy}

\def\O{\mathord{\mathbb O}}
\def\R{\mathord{\mathbb R}}
\def\H{\mathord{\mathbb H}}

\def\S{\mathord{\mathbb S}}
\def\3{{\ss}}
\def\2{\frac{1}{2}}
\def\x{\times}
\def\.{\cdot}

\def\<{\langle}
\def\>{\rangle}

\def\Aut{\mathop{\rm Aut}\nolimits}
\def\Im{\mathop{\rm Im\,}\nolimits}

\def\Ad{\mathop{\rm Ad}\nolimits}

\def\trace{\mathop{\rm trace\,}\nolimits}
\def\diag{\mathop{\rm diag}\nolimits}
\def\Spin{\mathop{\rm Spin}\nolimits}
\def\Fix{\mathop{\rm Fix}\nolimits}

\def\g{\mathfrak{g}}

\def\m{\mathfrak{m}}

\def\h{\mathfrak{h}}

\def\ms{\medskip\noindent}

\def\beq{\begin{equation}}
\def\eeq{\end{equation}}
\def\bea{\begin{eqnarray}}
\def\eea{\end{eqnarray}}

\def\esm{\end{smallmatrix}\right)}
\def\bsm{\left(\begin{smallmatrix}}
\def\epm{\end{pmatrix}}
\def\bpm{\begin{pmatrix}}

\author[J.-H.\ Eschenburg]{Jost-Heinrich Eschenburg}
\address{Institut f\"ur Mathematik, Universit\"at Augsburg, D-86135 Augsburg, Germany}
\email[Eschenburg]{eschenburg@math.uni-augsburg.de}

\author[T.\ Sakai]{Takashi Sakai}
\address{Department of Mathematical Sciences, Tokyo Metropolitan University, 1-1 Minami-Osawa, Hachioji, Tokyo 192-0397, Japan}
\email[Sakai]{sakai-t@tmu.ac.jp}

\title{Polars and antipodal sets for the outer $3$-symmetric space $\S^7\x\S^7$}

\thanks{The second author was supported by JSPS KAKENHI Grant Number JP21K03250.}



\begin{document}

\begin{abstract} 
We determine the polar and the maximal antipodal set $P$ for the outer 3-symmetric
space $\S^7 \x \S^7 = \Spin_8/G_2$ where the 3-symmetric structure is given by the
triality automorphism $\tau$ on $\Spin_8$. It turns out that $P$ has three elements.
The 3-symmetric structure extends to a (non-abelian) $S_3$-structure which we
also investigate.
\end{abstract}

\maketitle

\newtheorem{theorem}{\bf Theorem}[section]
\newtheorem{corollary}[theorem]{\bf Corollary}
\newtheorem{proposition}[theorem]{\bf Proposition}
\newtheorem{lemma}[theorem]{\bf Lemma}
\newtheorem{definition}[theorem]{\bf Definition}

\setcounter{tocdepth}{3}

\section{Introduction}

The Riemannian product $\S^{n-1}\x\S^{n-1}$ of two round $(n-1)$-spheres with equal size
is a symmetric space and allows several obvious $\Gamma$-symmetric structures. However,
for $n=8$ there is an outer 3-symmetric structure (even a non-abelian $S_3$-structure) 
which is not at all obvious. It is related to properties of the octonions, 
in particular to the triality automorphismus 
$\tau$ making the rotation of the Dynkin diagram $D_4$ for $\Spin_8$. According to
\cite[Thm.5.5,\,p.107]{WG} (see also \cite{SY}) this is one of the two outer 3-symmetric
spaces $G/H$ with $G$ simple, where $G = \Spin_8$ in both cases. We compute the polar
and the maximal antipodal set $P$ for this structure; $P$ consists
of three points and is preserved by the $S_3$-structure.

\section{$\Gamma$-symmetric spaces}
Let $G$ be a compact connected Lie group with left-invariant Riemannian metric
and $\Gamma\subset \Aut(G)$ a subgroup preserving the metric. Let $H\subset G$
be a closed subgroup with
$$
	\Fix(\Gamma)^e \subset H \subset \Fix(\Gamma)
$$
where $\Fix(\Gamma) = \{g\in G: \gamma(g) = g\ \forall\,{\gamma\in\Gamma}\}$ is the fixed group
of $\Gamma$ and $\Fix(\Gamma)^e$ its identity component. Suppose that the metric on $G$
is $H$-right-invariant.
Then $X := G/H$ inherits a $G$-invariant submersion metric from $G$. Further,
the action of $\Gamma$ on $G$ descends to an isometric action on $X$
by setting $\gamma(gH) := \gamma(g)H$. This is well defined 
since $\gamma(gh) = \gamma(g)h$ for all $h\in H\subset G^\Gamma$. Further,
$o = eH\in X$ is an isolated fixed point of this action. In fact, let $\h\subset\g$
be the Lie algebra of $H$ and $\m = \h^\perp$ the orthogonal complement 
which will be identified with the tangent space $T_oX$. Then $\Gamma$ acts on $\m$ 
without fixed vectors, hence $o$ is an isolated fixed point of $\Gamma$ on $X$,
and $X$ is called a generalized $s$-manifold (\cite{OS}), in particular,
if $\Gamma$ is finite abelian, then $X$ is called a $\Gamma$-symmetric space (cf.\ \cite{L}, \cite{BG}).

To generalize symmetric spaces we need a conjugate group $\Gamma_x\subset I(X)$ 
with isolated fixed point $x$ for every $x\in X$. In fact, $\Gamma_x = g\Gamma g^{-1}$
when $x = gH$. This makes sense since we consider both $G$ and $\Gamma$ 
as subgroups of the isometry group $I(X)$ where they generate a subgroup $G\.\Gamma\subset I(X)$.
This is a semidirect product of $G$ and $\Gamma$ 
since for all $g\in G$ and $\gamma\in\Gamma$ we have
\beq \label{gammag}
	\gamma g = g^\gamma\gamma
\eeq
where $g^\gamma := \gamma(g)$. Note that $(g^\gamma)^\delta = \delta(\gamma(g)) = g^{\delta\gamma}$.

More precisely, for any $x\in X$ we choose some $g_x\in G$ with $x = g_xH$ 
and define a group homomorphism
\beq \label{phix}
	\phi_x : \Gamma \to I(X): \gamma\mapsto \gamma_x:= g_x\gamma g_x^{-1} 
\eeq
This depends only on $x$: For all $h\in H$ we have $h\gamma h^{-1} = hh^{-1}\gamma = \gamma$,
and replacing $g_x$ by $\tilde g_x = g_xh$ we obtain $\tilde g_x\gamma\tilde g_x^{-1} =
g_xh\gamma h^{-1}g_x^{-1} = \gamma_x$. 

Using \eqref{gammag} we see the so called kai property
\beq \label{conj}
	\gamma_x\delta_y\gamma_x^{-1} = \phi_{\gamma_x(y)}(\gamma\delta\gamma^{-1})
\eeq
for all $x,y\in X$ and $\gamma,\delta\in\Gamma$ (see \cite{OS}).

\section{The triality automorphism}

The only compact simple group $G$ with $|\Aut(G)/\Aut(G)^e| > 2$ is $\Spin_8$. 
This is due to the triality automorphism which we want to explain next, following \cite{EQT}.

The Clifford algebra $Cl_8$ is the algebra of real $16\x 16$ matrices
$$Cl_8 = \left\{\bsm A & B \cr C & D \esm: A,B,C,D \in\R^{8\x 8}\right\}.$$
The even and odd parts $Cl_8^\pm$ consist of the block diagonal and antidiagonal
matrices $\bsm A &  \cr  & D \esm$ and $\bsm  & B \cr C &  \esm$.
The Euclidean vector space $\R^8$ can be viewed as the division algebra $\O$ of the octonions
where the basis vector $e_1 = (1,0,\dots,0)$ becomes the unit $1$. This
is embedded in $Cl_8^-$ by the linear isometry
\beq
	\R^8 = \O\ni x \mapsto  \hat x := \bpm 0 & -L(\overline x) \cr L(x) & 0 \epm
\eeq 
where $L(x)$ denotes the left translation on $\O$ with $x\in\O$, 
where $\kappa : x\mapsto\overline x$, $\kappa = \diag(1,-1,\dots,-1)$ is the octonionic
conjugation  (an anti-automorphism on $\O$) and where the inner
product on $\R^{16\x 16}$ is the normalized trace metric:
$$\begin{matrix}
	\<A,B\> = \frac{1}{16}\trace A^tB.
\end{matrix}$$
The group $\Spin_8$ is the subgroup of $SO_8 \x SO_8$ which preserves this
set of matrices under conjugation: $(A,B) \in \Spin_8$ $\iff$ 
$\Ad\bsm A\cr & B\esm \hat x = \hat w$ for some $w\in\O$ $\iff$
$$
	\bpm A & \cr & B \epm \bpm & -L(\overline x) \cr L(x) & \epm \bpm A^t & \cr & B^t \epm
	= \bpm & -L(\overline w) \cr L(w) & \epm
$$
for some $w \in \O$. This is equivalent to the equation 
\beq	\label{Bx=wA}
	B\.L(x) = L(w)\.A.
\eeq
From $BL(x)A^t = L(w)$ it follows that the map $C: x\mapsto w$ is linear and isometric, 
hence $C\in SO_8$ (using connectivity). 
Applying both sides of \eqref{Bx=wA} to some $y\in\O$ we see that $(A,B) \in \Spin_8$ if and only if there is some $C\in SO_8$ with 
\beq \label{ABC}
  B(x\.y) = (Cx)\.(Ay)  
\eeq
for all $x,y \in \O$; this property is called {\em triality}. 

We may view $\Spin_8$ as the set of triples
$(A,B,C) \in (SO_8)^3$ with (\ref{ABC}).

\begin{lemma}
These triples form a subgroup of $(SO_8)^3$.
\end{lemma}

\proof
We show: When \eqref{ABC} holds for $(A,B,C)$ and $(A',B',C')$, then it also holds
for the product triple $(A'A,B'B,C'C)$. In fact, 
$B'B(x\.y) = B'(Cx\.Ay) = C'Cx\.A'Ay$.
Further, when $(A,B,C)\in \Spin_8$, the same is true for its inverse $(A^t,B^t,C^t)$,
that is $B^t(xy) = C^tx\.A^ty$. Indeed, 
$B(C^tx\.A^ty) = CC^tx\.AA^ty = xy = BB^t(xy)$.
\endproof

The {\em triality automorphism} of $\Spin_8$ is essentially a cyclic change in every triple. 
More precisely: 

\begin{lemma} 
Let $\kappa = \diag(1,-1,...,-1)$, $\kappa(x) = \overline x$ be the octonionic conjugation on $\O$.
Then $\tilde D := \kappa D \kappa \in SO_8$ for all $D\in SO_8$, and 
\beq \label{tau}
\tau : (A,B,C) \mapsto (\tilde B,\tilde C,A)
\eeq
is an automorphism of $\Spin_8$ with order 3.
\end{lemma}

\proof
Let $x,y\in\O$ with $|x| = |y| = 1$ and $z = xy$, thus
$|z| = 1$. Using quaternionic conjugation 
(an anti-automorphism of $\O$), we obtain
$$
(Cx)(Ay) = Bz \iff (\overline{Ay})(\overline{Cx}) = \overline{Bz}
	      \iff \overline{Cx} = (Ay)(\overline{Bz}).
$$
Now  $\overline{Cx} = \tilde C \overline x$  with  $\tilde C = \kappa C\kappa$.  
Thus (\ref{ABC}) $\iff (Ay)(\tilde B\overline z) = \tilde C\overline x$,
and since $y\overline z = \overline x$ when $z,x$ have unit length, we have 
$(\tilde B,\tilde C,A) \in \Spin_8$.

Clearly, permutations and the maps $B\mapsto\tilde B$ etc.\ are automorphisms.

The iteration of $\tau$ gives
\beq \label{tau2}
(A,B,C) \mapsto (\tilde B,\tilde C,A) \mapsto (C,\tilde A,\tilde B) 
\mapsto (A,B,C). 
\eeq
\endproof

\begin{lemma}
There is an outer involution $\sigma$ on $\Spin_8\subset (SO_8)^3$,
\beq	\label{sigma}
	\sigma : (A,B,C) \mapsto (B,A,\tilde C).
\eeq
\end{lemma}

\proof Instead of
$Cx\.Ay = B(xy)$ we have to show $\tilde Cx\,B\tilde y = A(x\tilde y)$ for all $x,\tilde y\in\O$.
We transform the first equation as follows. Replacing $x$ by $\overline x$
and assuming that $x,y$ are unit vectors we obtain
\bea
C\overline x\,Ay = B(\overline xy) &\iff& C\overline x = B(\overline xy)\,\overline{Ay}	\cr
&\iff& \overline{C\overline x} = Ay\,\overline{B(\overline xy)}	\cr
&\iff& \tilde C(x)\.B(\overline xy) = Ay 	\nonumber
\eea
Setting $\tilde y = \overline xy$ we have $x\tilde y = x\overline xy = y$ and
hence $\tilde C x\,B\tilde y = A(x\tilde y)$ as required.
\endproof

\begin{lemma}
The two automorphisms $\sigma$ and $\tau$ generate a group $\hat\Gamma\subset\Aut(\Spin_8)$
which is isomorphic to $S_3$.
\end{lemma}

\proof
We have to check the relation $\sigma\tau\sigma = \tau^{-1}$. In fact,
$$
(A,B,C) \buildrel \sigma\over\mapsto (B,A,\tilde C) 
\buildrel\tau\over\mapsto (\tilde A,C,B) 
\buildrel\sigma\over\mapsto (C,\tilde A,\tilde B),
$$
and $(C,\tilde A,\tilde B)$ is the triple for $\tau^2 = \tau^{-1}$ in \eqref{tau2}.
\endproof

\begin{lemma}
The fixed group of $\hat\Gamma = \<\tau,\sigma\>_{\rm grp}$ and of $\Gamma = \<\tau\>_{\rm grp}$
acting on $\Spin_8$ is $G_2 = \Aut(\O)$.
\end{lemma}

\proof
Every triple $(A,B,C)\in \Spin_8$ is fixed by $\tau$ or by both $\tau$ and $\sigma$ $\iff$ $A=B=C$
$\iff$ $A \in \Aut(\O)$, using \eqref{ABC}.
\endproof

\section{The outer $\hat\Gamma$-symmetric space $\S^7\x\S^7$.}

As a consequence of the last sections, $X = \Spin_8/G_2$ is an outer $\hat\Gamma$-symmetric space.
The group $\Spin_8$ can be viewed as a subgroup of $SO_8\x SO_8$ since $C$ is determined
by $A,B$ via \eqref{ABC}. As such it acts naturally on $\S^7\x\S^7 \subset\O^2$.

\begin{lemma} The isotropy group of this action at the point 
$o := (1,1)\in \S^7\x\S^7=X$ is $G_2 = \Aut(\O)$. Consequently, $X \cong \S^7\x\S^7$ 
$(\Spin_8$-equivariantly diffeomorphic$)$.
\end{lemma}

\proof
When $(A,B,C)\in \Spin_8$ with $A(1) = B(1) = 1$, then from \eqref{ABC}
we obtain $C(1)\.A(1) = B(1)$, hence $C(1) = 1$, further $C(1)\.Ay = By$ whence $A=B$,
and $Cx\.A(1) = By$, thus $C = B$. Therefore $(A,B,C) = (A,A,A)$ 
which implies $A \in G_2$ by \eqref{ABC}.\\ 
Since $\dim(\Spin_8/G_2) = 28-14 = 14$
is the dimension of $\S^7\x\S^7$, the action of $\Spin_8$ on $\S^7\x\S^7$ is transitive. 
The $\Spin_8$-equivariant diffeomorphism is
\beq	\label{S7S7}
	\Spin_8/G_2 \ni (A,B)\.G_2 \mapsto (A(1),B(1))\in \S^7\x\S^7.
\eeq
\endproof

How does the action of $\hat\Gamma$ on $\Spin_8/G_2$ look like on $\S^7\x\S^7$?
We have $\tau(A(1),B(1)) = (\tilde B(1),\tilde C(1))$ by \eqref{tau}. From \eqref{ABC} we have
$C(1)\.A(1) = B(1)$, thus $C(1) = B(1)\.\overline{A(1)}$. Setting $x = A(1)$ and $y = B(1)$
we have $\tau(x,y) = (\overline y,\overline{y\.\overline{x}})$, hence
\beq \label{tauS}
	\tau(x,y) = (\overline y,x\overline y).
\eeq
Once more we can see the order 3 of this map:
$$
(x,y)
\mapsto (\overline y,x\overline y)
\mapsto (\overline{x\overline y},\overline y\, \overline{x\overline y}) = (y\overline x,\overline x)
\mapsto (x,y\overline x x) = (x,y).
$$

\subsection{The fixed points of $\hat\Gamma$} \label{fixed}
Thus any $(x,y)\in \S^7\x\S^7$ is fixed by $\tau$ $\iff$
$x=\overline y$ and $y = x\overline y$ $\iff$ $\overline x = y = x^2$.
In particular, $x^2 = \overline x$, thus $x^3 = 1$, and $y = \overline x$. Hence
the fixed set of $\tau$ in $X = \S^7\x\S^7$ is 
\beq \label{Fixtau}
	\Fix(\tau) = \{(x,\overline x):x\in\O,\,x^3=1\} = \{(1,1)\} \cup Y,
\eeq
$$
\begin{matrix}
	Y = \{(x,\overline x): x = \2(-1+\sqrt3v): v\in\S^6\}
\end{matrix}
$$
where $\S^6$ is the unit sphere in $\R^7 = \Im(\O) = \O\ominus\R$.\goodbreak

Further, $\sigma$ acts on $X$ by $(A(1),B(1))\mapsto (B(1),A(1))$, cf.\ \eqref{sigma}, 
hence its fixed set is the diagonal,
$$
	\Fix(\sigma) = \{(x,x): x\in\S^7\}.
$$ 
Since $\Fix(\sigma)\cap Y = \emptyset$, we obtain
\beq	\label{FixGamma}
	\Fix(\hat\Gamma) = \Fix(\sigma)\cap \Fix(\tau) = \{(1,1)\}.
\eeq

\section{Maximal antipodal sets for $\Gamma = \<\tau\>$}

Let $X$ be a $\Gamma$-symmetric space. An {\em antipodal set} in $X$ is a subset $P\subset X$
such that $P \subset \Fix(\Gamma_p)$ for all $p\in P$. By translating $P$ if necessary 
we may always assume $o\in P$. Clearly for the $S_3$-symmetric space $\S^7\x\S^7$, the
only antipodal set is $\{o\}$. In the following we will consider $\S^7\x\S^7$ as a
$\Gamma$-symmetric space for $\Gamma = \<\tau\>_{\rm grp}$, an outer 3-symmetric space.
We have already seen that $Y$ is a {\em polar} for $o$, a positive-dimensional
component of $\Fix(\Gamma)$ (mind that $\Gamma = \Gamma_o$).

First let us recall a well known facts on octonions.

\begin{lemma}
For all $s,x\in\O$ we have
\beq \label{sxs}
	L(s)L(x)L(s) = L(sxs).
\eeq
\end{lemma}

\proof
Since any two octonions lie in a commom quaternionic subalgabra, we may assume $s,x\in\H\subset\O$.
Clearly, \eqref{sxs} holds on $\H$. We have to show it on $\ell\H$ where $i,j,\ell$ are the usual
unit generators of $\O$ (with $\ell\perp i,j,ij$). 
Recall for every $h\in\H$:
$$
	L(h)L(\ell) = L(\ell) L(\overline h).
$$
In fact, let $h = h_o + h'$ with $h_o\in\R$ and $h'\perp \R$. Then
$$h_o\ell+h'\ell = \ell h_o -\ell h' = \ell \overline h.$$
Therefore we have on $\H$:
$$
L(s)L(x)L(s)L(\ell) = L(\ell)L({\overline s})L({\overline x})L({\overline s}) 
= L(\ell) L(\overline s\overline x\overline s) = L({sxs})L(\ell).
\eqno{\qed}
$$

\begin{theorem}
Let $\Gamma\subset\Aut(\Spin_8)$ be the group generated by the triality automorphism $\tau$.
In the outer $\Gamma$-symmetric space $\S^7\x\S^7 = \Spin_8/G_2$, 
every maximal antipodal set containing $o = (1,1)$ equals 
	$$P = \{(1,1),(s,\overline s),(\overline s,s)\}$$
with $s = \2(-1+\sqrt3\,v)$ for some arbitrary $v\in \S^6\subset\Im\O$. The outer involution
$\sigma$ fixes $(1,1)$ and interchanges the other two elements of $P$.
\end{theorem}

\proof
Any $p=(s,\overline s)\in \Fix(\tau)\subset \S^7\x\S^7$ 
can be written in the form $go$ with $o=(1,1)$ and $g = (A,B)$ for 
	$$A = L(s), \ \ \ B = L(\overline s).$$ 
Then
$g = (A,B) \in \Spin_8$ since \eqref{ABC} holds for $A,B,C$ with $Cx = \overline sx\overline s$: 
\bea
(Cx)(Ay)
&=& (\overline s x \overline s)(sy)
= L({\overline s x \overline s})L(s)y \cr
&\buildrel \eqref{sxs}\over =& L({\overline s})L(x)L({\overline s})L(s)y
= L({\overline s})L(x)y
= B(xy).
\nonumber\eea
Clearly, both $\Gamma$ and $\Gamma_p = g\Gamma g^{-1}$ are fixing $o$ and $p$,
and $\Fix(\Gamma_p) = g\Fix(\Gamma)$ contains $p = go$ and $o$. Note that
$o = g(q)$ for $q = (\overline s,s)$, and $(\overline s,s) \in Y \subset \Fix(\Gamma)$,
see \eqref{Fixtau}. Therefore $q\in\Fix(\Gamma)\cap \Fix(\Gamma_p)$. Applying $\sigma$
we see also $p\in \Fix(\Gamma)\cap\Fix(\Gamma_q)$, and $\Gamma$ fixes $p,q\in Y$.
Therefore $P$ is antipodal.

A {\em maximal} antipodal set through $o$ and $p$ might contain other points
$\tilde q\in \Fix(\Gamma)\cap g\Fix(\Gamma)$. 
The elements of $g\Fix(\Gamma)$ with $g = (L(s),L(\overline s))$ are 
$\tilde q = (st,\overline s\overline t)$
for any $t\in\O$ with $t^3 = 1$, see Subsection \ref{fixed}.
Such $\tilde{q}=(st,\overline s\overline t)$ lies in $\Fix(\Gamma)$ if and only if 
$\overline s\overline t = \overline{st}$ 
which means that $s$ and $t$ commute and hence $(st)^3 = s^3t^3 = 1$. 
There are three solutions: $t \in \{1,s,\overline s\}$, 
leading to $\tilde q = (u,\overline u)$ for $u \in \{s,\overline s,1\}$. 
The only new solution is $u = \overline s$ with $\tilde q = (\overline s,s) = q$.
\endproof

\ms \hskip 4cm
\includegraphics{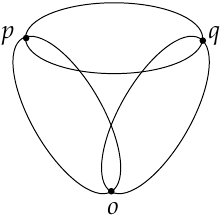}

\ms
The fixed space of $\Gamma = \Gamma_o$ consists of the isolated point $o$ 
and the 6-sphere $Y$ opposite to $o$, a polar for $o$, 
drawn as an ellipse in the preceding figure. The other points $p$ and $q$ in $P$
are antipodal in $Y$. Likewise, the fixed sets of $\Gamma_p$
and $\Gamma_q$ consist $p$ or $q$ (resp.) together with opposite 6-spheres.
The three 6-spheres intersect pairwise in antipodal points. These three points $o,p,q$
form a maximal (in fact ``great'') antipodal set $P$ for the $\Gamma$-structure on $\S^7\x\S^7$.

\end{document}